\documentclass[12pt]{amsart}

\usepackage{graphicx}

\newtheorem{tw}{Theorem}
\newtheorem{Lem}[tw]{Lemma}

\newcommand{\cal}[1]{\mathcal{#1}}
\newcommand{\sig}{\sigma}
\newcommand{\eps}{\varepsilon}
\newcommand{\ro}{\varrho}
\newcommand{\bet}{\beta}

\newcommand{\te}{\theta}

\newcommand{\kre}[1]{\overline{#1}}
\newcommand{\gen}[1]{\langle #1 \rangle}
\newcommand{\map}[3]{#1\colon #2\to #3}
\newcommand{\Map}[2]{#1\colon #2\to #2}
\newcommand{\field}[1]{\mathbb{#1}}
\newcommand{\zz}{\field{Z}}
\newcommand{\rr}{\field{R}}

\newcommand{\lst}[2]{{#1}_1,\dotsc,{#1}_{#2}}
\newcommand{\Mpm}{{\cal{M}}_{g}^{\pm}}
\newcommand{\Mg}{{\cal{M}}_{g}}
\newcommand{\Mtw}{{\cal{M}}_{2}}
\newcommand{\Mtwpm}{{\cal{M}}_{2}^{\pm}}
\newcommand{\Mh}{{\cal{M}}^{h}_{g}}
\newcommand{\Mhpm}{{\cal{M}}^{h\pm}_{g}}

\newenvironment{dow}{\begin{proof}}{\end{proof}}

\begin{document}




\title{Small torsion generating sets for hyperelliptic mapping class groups}
\thanks{Supported by KBN 1 P03A 026 26}

\author{Micha\l\ Stukow}
\email{trojkat@math.univ.gda.pl}

\address{Institute of Mathematics, University of Gda\'nsk, Wita Stwosza 57,
80-952 Gda\'nsk, Poland }

\begin{abstract}
We prove that both the hyperelliptic mapping class group and the extended hyperelliptic
mapping class group are generated by two torsion elements. We also compute the index of
the subgroup of the hyperelliptic mapping class group which is generated by involutions
and we prove that the extended hyperelliptic mapping class group is generated by three
orientation reversing involutions.
\end{abstract}

%

\maketitle
\section{Introduction}
Let $S_g$ be a closed orientable surface of genus $g\geq 2$. Denote by $\Mpm$ the
\emph{extended mapping class group} i.e. the group of isotopy classes of homeomorphisms
of $S_g$. By $\Mg$ we denote the \emph{mapping class group} i.e. the subgroup of $\Mpm$
consisting of orientation preserving maps. We will make no distinction between a map and
its isotopy class, so in particular by the order of a homeomorphism $\Map{h}{S_g}$ we
mean the order of its class in $\Mpm$.

\begin{figure}[h]
\begin{center}
\includegraphics[width=0.75\textwidth]{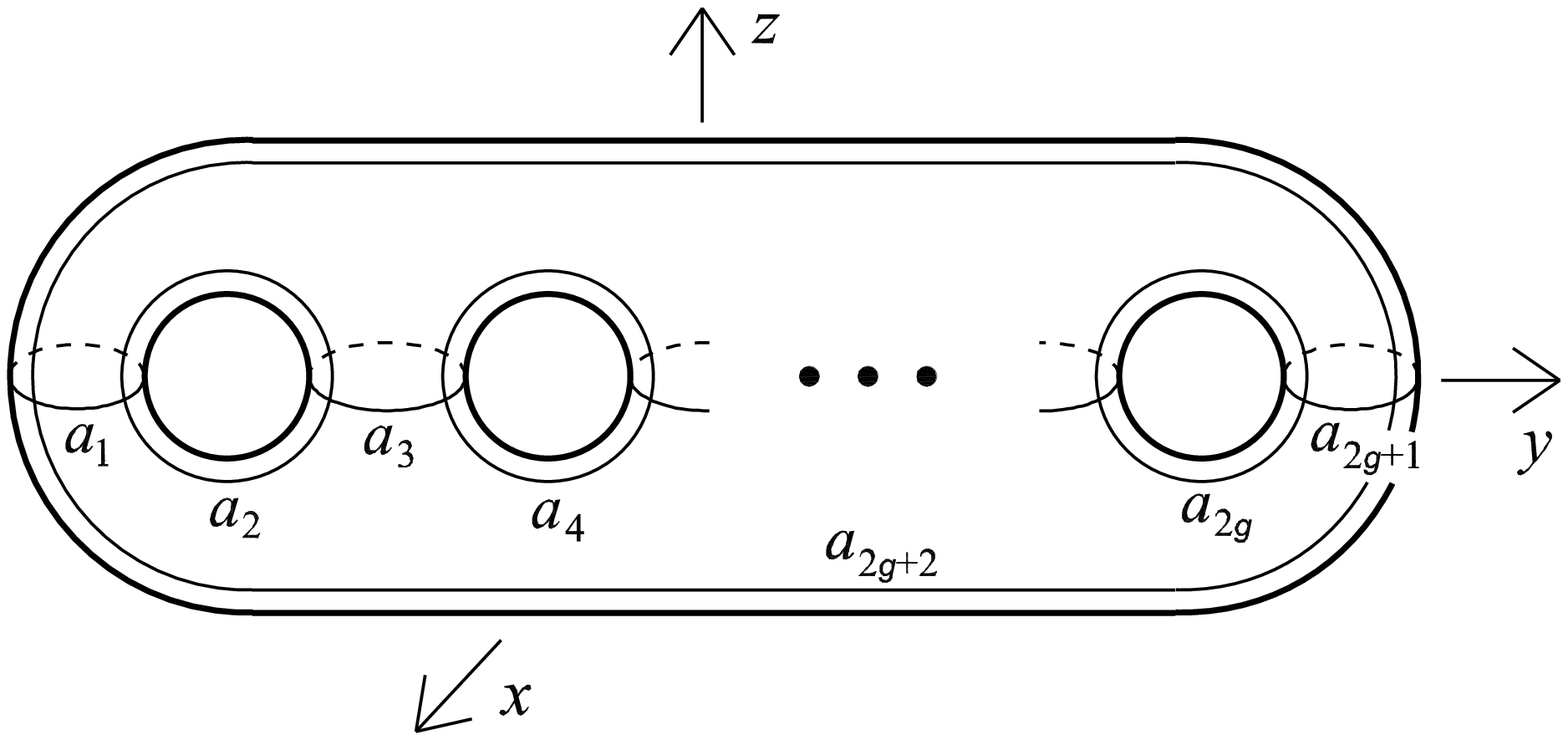}
\caption{Surface $S_g$ embedded in $\rr^3$.}\label{r1} %
\end{center}
\end{figure}

Suppose that $S_g$ is embedded in $\rr^3$ as shown in Figure \ref{r1}, in such a way that
it is invariant under reflections across $xy,yz,xz$ planes. Let $\Map{\ro}{S_g}$ be a
\emph{hyperelliptic involution}, i.e. the half turn about the $y$-axis.

The \emph{hyperelliptic mapping class group} $\Mh$ is defined to be the centraliser of
$\ro$ in $\Mg$. In a similar way, we define the \emph{extended hyperelliptic mapping
class group} $\Mhpm$ to be the centraliser of $\ro$ in $\Mpm$. For $g=2$ it is known that
$\Mtw={{\cal{M}}^{h}_2}$ and $\Mtwpm={{\cal{M}}^{h\pm}_2}$.

The problem of finding certain generating sets for groups $\Mg$ and $\Mpm$ has been studied
extensively -- see \cite{Bre-Farb,GG-MS,Hump,Kassabov,Kork1,MacMod,McCarPap1,inv,Wajn} and
references there. In particular it is known that the group $\Mg$ could be generated by $2$ torsion
elements \cite{Kork1}, and by small number of involutions \cite{Kassabov}. The problem of finding
the minimal generating set consisting of involutions is still open.

Similar results hold for the extended mapping class group, namely the group $\Mpm$ is generated by
$2$ elements \cite{Kork1}, and by $3$ symmetries (i.e. orientation reversing involutions)
\cite{inv}. The question if it is possible to generate $\Mpm$ by two torsion elements is open.

The purpose of this paper is to give a full answer to analogous questions in the case of the
hyperelliptic mapping class group and the extended hyperelliptic mapping class group, i.e. we will
show that both groups are generated by two torsion elements and $\Mhpm$ is generated by three
symmetries. It will be also observed that the subgroup $I_g$ of $\Mh$ which is generated by
involutions is a proper subgroup and we will compute the index $[\Mh:I_g]$.

The importance of torsion elements in the group $\Mg$ ($\Mpm$) follows from the fact that
any such element could be realized as an analytic (dianalytic) automorphism of some
Riemann surface. Similarly torsion elements in $\Mh$ ($\Mhpm$) correspond to analytic
(dianalytic) automorphisms of hyperelliptic Riemann surfaces, i.e. complex algebraic
curves with affine part defined by an equation $y^2=f(x)$, where $f$ is a polynomial with
distinct roots.

Notice that since $\Mtw={{\cal{M}}^{h}_2}$ and $\Mtwpm={{\cal{M}}^{h\pm}_2}$, our results
imply that $\Mtw$ and $\Mtwpm$ are generated by two torsion elements.
\section{Preliminaries}
Let $\lst{A}{2g+2}$ be the right Dehn twists along the curves $\lst{a}{2g+2}$ indicated in Figure
\ref{r1}. Denote also $B=A_1A_2\cdots A_{2g+1}$ and $\kre{B}=A_{2g+1}\cdots A_2 A_1$. It is known
\cite{Bir-Hil} that
$$\ro=B\kre{B}=A_1A_2\cdots A_{2g+1} A_{2g+1}\cdots A_2 A_1$$
and $\Mh$ admits the presentation:\\
\emph{generators:} $\lst{A}{2g+2},B,\ro$\\
\emph{defining relations:}
\begin{align}
 &\ro=A_1A_2\cdots A_{2g+1} A_{2g+1}\cdots A_2 A_1 \label{roro}\\
 &B=A_1A_2\cdots A_{2g+1}\label{BB}\\
 &A_{j}=BA_{i}B^{-1} \quad \text{$j\equiv {i+1}\mod {2g+2} $}\label{B}\\
 &A_iA_j=A_jA_i, \quad \text{$2\leq |i-j|\leq 2g$}\\
 &A_iA_jA_i=A_jA_iA_j \quad \text{$j\equiv {i+1}\mod {2g+2} $}\label{braid}\\
 &B^{2g+2}=1\label{rzad}\\
 &\ro^2=1\\
 &\ro A_i=A_i \ro\label{com}
\end{align}
The above presentation follows from the presentation in \cite{Bir-Hil} by adding generators
$A_{2g+2},B,\ro $ and some superfluous relations. The relation \eqref{B} follows from easily
verified observation that
$$B(a_i)=a_j \quad \text{for $j\equiv {i+1}\mod {2g+2}$}$$
Combining this with relation \eqref{rzad} we have that $B$ has order $2g+2$. Observe also, that for
$k\in\zz$
$$B=B^k B B^{-k}=B^k A_1A_2\cdots A_{2g+1} B^{-k}=A_{1+k}A_{2+k}\cdots A_{2g+1+k}$$
where subscripts should be reduced modulo $2g+2$.

Let us also point out that the relation \eqref{B} implies that $\gen{B,A_i}=\Mh$ for every $1\leq
i\leq 2g+2$.

\section{Minimal torsion generating sets for $\Mh$ and $\Mhpm$}
\begin{tw}
For every $g\geq 2$ the group $\Mh$ is generated by two elements of order $2g+2$ and $4g+2$
respectively.
\end{tw}
\begin{dow}
Let $M=A_2A_3\cdots A_{2g+1}$. It is well known that $M$ has order $4g+2$ (cf.
\cite{Bre-Farb,inv}). Since $A_{1}=BM^{-1}$, we have $$\gen{B,M}=\gen{B,A_{1}}=\Mh$$
\end{dow}
Let $\sig$ be the reflection across the $yz$-plane (Figure \ref{r1}). Since $\sig(a_i)=a_i^{\pm
1}$, $1\leq i\leq 2g+2$ and $\sig$ reverses orientation, we have
\begin{equation}\label{sig}\sig
A_i\sig=A_i^{-1}\quad \text{for $1\leq i\leq 2g+2$}\end{equation} Therefore
\begin{equation}\sig B\sig=\kre{B}^{-1}=\ro B
\label{sigB}\end{equation}
\begin{Lem}
The order of the element $\bet=\sig B$ is finite and equal to $2g+2$ for $g$ odd and $4g+4$ for $g$
even.
\end{Lem}
\begin{dow}
Since $\bet^2=\sig B\sig B=\ro B^2$, we have $\bet^{2g+2}=\ro^{g+1}$, which completes the proof.
\end{dow}
\begin{Lem}
The order of the element $N=\sig A_{2g+1}^{-1}A_1A_2A_1^{-1}BA_{2g+1}^{-1}$ is finite and equal to
$2g$ for $g$ odd and $4g$ for $g$ even.
\end{Lem}
\begin{proof}
Using relations \eqref{B}--\eqref{sigB}, we compute
\begin{equation}\begin{split}N^2&=(\sig A_{2g+1}^{-1}A_1A_2A_1^{-1}BA_{2g+1}^{-1})
(\sig A_{2g+1}^{-1}A_1A_2A_1^{-1}BA_{2g+1}^{-1})\\
&=A_{2g+1}A_1^{-1}A_2^{-1}A_1\ro B
A_1A_2A_1^{-1}BA_{2g+1}^{-1}\\
&=A_{2g+1}(A_1^{-1}A_2^{-1}A_1)
A_2A_3A_2^{-1}B^2A_{2g+1}^{-1}\ro\\
&=A_{2g+1}(A_2A_1^{-1}A_2^{-1})
A_2A_3A_2^{-1}B^2A_{2g+1}^{-1}\ro\\
&=A_{2g+1}(A_2A_3)(A_1^{-1}A_2^{-1})B^2A_{2g+1}^{-1}\ro
\end{split}\label{S2}
\end{equation}
Similar computations show that
\begin{equation*}
\begin{split}
N^{2g}&=A_{2g+1}(A_2\cdots A_{2g}A_{2g+1})(A_1^{-1}\cdots
A_{2g-1}^{-1}A_{2g}^{-1})B^{2g}A_{2g+1}^{-1}\ro^{g}\\
&=A_{2g+1}(A_2\cdots A_{2g}A_{2g+1}A_{2g+2})(A_{2g+2}^{-1}A_1^{-1}\cdots
A_{2g-1}^{-1}A_{2g}^{-1})B^{2g}A_{2g+1}^{-1}\ro^{g}\\
&=A_{2g+1}B\ro B B^{2g}A_{2g+1}^{-1}\ro^{g}=\ro^{g+1}\end{split}\end{equation*} This completes the
proof.
\end{proof}
Let $G=\gen{\bet,N}$. We will show that $G=\Mhpm$, but first we need the following lemma.
\begin{Lem}\label{Lem}
If $g\geq 3$ then $A_{2g+1}A_1^{-1}\in G$.
\end{Lem}
\begin{dow}
Suppose first that $g\geq 4$. Using \eqref{S2} we have
\begin{equation*}
\begin{split}
&N^{-2}(\bet^4N^2\bet^{-4})(\bet^{-2}N^2\bet^2)(\bet^2 N^{-2} \bet^{-2})\\
=&(A_{2g+1}B^{-2}A_2 A_1 A_3^{-1} A_2^{-1} A_{2g+1}^{-1})\ro(A_{3}A_6 A_7 A_5^{-1}
A_6^{-1} B^2 A_{3}^{-1})\ro\\ & (A_{2g-1}A_{2g+2} A_1 A_{2g+1}^{-1} A_{2g+2}^{-1} B^2
A_{2g-1}^{-1})\ro
(A_{1}B^{-2}A_4 A_3 A_5^{-1} A_4^{-1} A_{1}^{-1})\ro\\
=&A_{2g+1}A_{2g+2} A_{2g+1} A_1^{-1} A_{2g+2}^{-1} A_{2g-1}^{-1}A_{1}(A_4 A_5 A_3^{-1} A_4^{-1}
A_{3}^{-1}) \\  &(A_{2g-1}A_{2g+2} A_1 A_{2g+1}^{-1} A_{2g+2}^{-1} A_{2g+1}^{-1})A_{3}A_4
A_3A_5^{-1}A_4^{-1} A_{1}^{-1}\\
=&A_{2g+1}A_{2g+2} A_{2g+1} A_1^{-1} A_{2g+2}^{-1} A_{2g-1}^{-1}A_{1}(A_{2g-1}A_{2g+2} A_1
A_{2g+1}^{-1} A_{2g+2}^{-1} A_{2g+1}^{-1}) \\  & (A_4 A_5 A_3^{-1} A_4^{-1} A_{3}^{-1})A_{3}A_4 A_3
A_5^{-1}A_4^{-1} A_{1}^{-1}\\
=&A_{2g+1}A_{2g+2} A_{2g+1} A_1^{-1} A_{2g+2}^{-1} (A_{1}A_{2g+2} A_1) (A_{2g+1}^{-1} A_{2g+2}^{-1}
A_{2g+1}^{-1})A_{1}^{-1}\\
=&A_{2g+1}A_{2g+2} A_{2g+1} A_1^{-1} A_{2g+2}^{-1} (A_{2g+2}A_{1} A_{2g+2}) (A_{2g+2}^{-1}
A_{2g+1}^{-1} A_{2g+2}^{-1})A_{1}^{-1}\\
=&A_{2g+1}A_1^{-1}
\end{split}
\end{equation*}
If $g=3$ similar, but rather long computations\footnote{There are available at the URL:\\
\texttt{http://www.math.univ.gda.pl/\~{}trojkat/comgenhi.pdf}} show that
$$\bet^{-4}N^{-2}\bet N^{-2}\bet^{-1}N^{-1}\bet^2 N^2\bet^{-4} N\bet N^{-3}\bet^4 N^{-1}\bet=A_7A_1^{-1}$$
\end{dow}
\begin{tw}
For every $g\geq 2$ the group $\Mhpm$ is generated by two elements of finite order.
\end{tw}
\begin{dow}
Let $\bet$ and $N$ be elements defined above. Since $\bet$ satisfies the relation
$$A_{j}^{-1}=\bet A_{i} \bet^{-1}\quad \text{$j\equiv {i+1}\mod {2g+2} $}$$
to prove the theorem it is enough to show that $A_i\in G=\gen{\bet, N}$ for some $1\leq i\leq
2g+2$.

If $g\geq 3$ then from Lemma \ref{Lem} and the above relation follows that
$$A_{2g+1}A_{1}^{-1},A_{2g}A_{2g+2}^{-1},A_{2g-1}A_{2g+1}^{-1}\in G$$
Therefore
\[\begin{split}G\ni &(A_{2g-1}A_{2g+1}^{-1})(A_{2g+1}A_1^{-1})(A_{2g+2}^{-1}A_{2g})\bet^{-1}N(A_{2g+1}A_{2g-1}^{-1})\\
&=A_{2g-1}A_1^{-1}A_{2g+2}^{-1}A_{2g}B^{-1}A_{2g+1}^{-1}A_1A_2A_1^{-1}BA_{2g+1}^{-1}A_{2g+1}A_{2g-1}^{-1}\\
&=A_{2g-1}A_1^{-1}A_{2g+2}^{-1}A_{2g}A_{2g}^{-1}A_{2g+2}A_1A_{2g+2}^{-1}A_{2g-1}^{-1}\\
&=A_{2g-1}A_{2g+2}^{-1}A_{2g-1}^{-1}=A_{2g+2}^{-1}\\
\end{split}\]
If $g=2$ then one could verify that
$$N^{-1}\bet N^2\bet N\bet N^{-1}\bet^{-1}N^{-1}\bet N^2\bet N\bet=A_3^{-1}$$
\end{dow}
\section{Involutions as generators for $\Mh$ and $\Mhpm$}
Observe that from a presentation for the group $\Mh$ we obtain
$$H_1(\Mh,\zz)=\begin{cases}\zz_{4g+2}&\text{for $g$ even}\\
\zz_{8g+4}&\text{for $g$ odd}
\end{cases}$$
Since neither $\zz_{4g+2}$ nor $\zz_{8g+4}$ is generated by involutions the same
conclusion holds for $\Mh$.

Denote by $I_g\leq \Mh$ the subgroup generated by involutions. Clearly this is a normal
subgroup of $\Mh$. Our next goal is to describe the quotient $\Mh/I_g$. To achieve it we
will follow similar lines to \cite{McCarPap1}.

Let $S$ be the half turn $S$ about the $z$-axis (Figure \ref{r1}).

\begin{Lem}\label{ConInv}
The number of conjugacy classes of involutions in $\Mh$ is equal to $2$ for $g$ even and $3$ for
$g$ odd. These classes are represented by $\ro,S$ and $\ro,S,\ro S$ respectively.
\end{Lem}
\begin{dow}
Since $\ro$ is central, its conjugacy class consists of one element. Therefore we could restrict
ourselves to conjugacy classes of involutions different from $\ro$.

Let $R\in \Mh$ be an involution and $H=\gen{R,\ro}$. In particular $H$ has order $4$ and contains
$\ro$. From Theorem $4$ in \cite{MaxHyp} follows that there are exactly $2$ conjugacy classes of
subgroups of $\Mh$ having these two properties. To identify these conjugacy classes, it is enough
to find two nonconjugate examples of such subgroups. The first one is the dihedral group
$\gen{S,\ro}$. An example of a subgroup in the second conjugacy class follows from easily verified
fact that
$$(A_1A_{2g+2}^{-1}B)^{2g}=\ro$$
In particular the group $\gen{(A_1A_{2g+2}^{-1}B)^g}$ is a cyclic group of order $4$ representing
the second conjugacy class. Therefore $H$, as a dihedral group, is conjugate to $\gen{S,\ro}$,
hence $R$ is conjugate to either $S$ or $\ro S$. To complete the proof it is enough to show that
$S$ and $\ro S$ are conjugate if and only if $g$ is even.

First suppose that $g$ is even. Let $t$ be a circle fixed by $S$, and
$$\te=(A_{g+2}A_{g+3}\cdots A_{2g}A_{2g+1})^{g+1}$$
From geometric point of view, $\te$ is a half-twist about $t$, i.e. it is a half-turn of the right
half of $T_g$, in particular $\te^2=T$ -- twist along $t$. From this geometric interpretation it is
clear that $$\te(S\te^{-1} S^{-1})=\ro$$ Hence $\te S\te^{-1}=\ro S$.

On the other hand, if $g$ is odd, then automorphisms induced by $S$ and $\ro S$ on
$H_1(T_g,\zz)$ have different eigenvalues, so $S$ and $\ro S$ can not be conjugate in
$\Mh$.
\end{dow}
\begin{Lem}\label{cyc}
The quotient $\Mh/I_g$ is cyclic.
\end{Lem}
\begin{dow}
Since $S(a_1)=a_{2g+1}$ we have
$$A_1A_{2g+1}^{-1}=A_1(SA_{1}^{-1}S^{-1})=(A_1SA_1^{-1})S^{-1}\in I_g$$

Now observe that for any $3\leq i\leq 2g$, we can construct an element $F\in\Mh$ such that
$F(a_1)=a_1$ and $F(a_i)=a_{2g+1}$. In fact, if we define $F_j=A_jA_{j+1}$ for $3\leq j\leq 2g$
then $F_j(a_1)=a_1$, $F_j(a_j)=a_{j+1}$ so we could take $F=F_{2g}F_{2g-1}\cdots F_{i}$. Hence
$$A_1A_{i}^{-1}=F^{-1}(A_1A_{2g+1}^{-1})F\in I_g\quad\text{for $3\leq i\leq 2g$}$$
Finally we have
$$A_{2g+1}A_2^{-1}=S(A_1A_{2g}^{-1})S^{-1}\in I_g $$
Therefore all twists $\lst{A}{2g+1}$ are equal modulo $I_g$. Since they generate $\Mh$ this implies
that $\Mh/I_g$ is cyclic.
\end{dow}
\begin{tw}
The index $[\Mh:I_g]$ of the subgroup generated by involutions is equal to $2g+1$ for $g$
even and $4g+2$ for $g$ odd.
\end{tw}
\begin{dow}
Let $\map{\pi}{\Mh}{H_1(\Mh,\zz)}$ be the canonical projection. By Lemma \ref{cyc}, $[\Mh,\Mh]\leq
I_g$, so
\begin{equation}\label{eq:pi}
[\Mh:I_g]=[H_1(\Mh,\zz):\pi(I_g)]\end{equation} From the presentation for $\Mh$, we have that for
any $1\leq i\leq 2g+2$ the group $H_1(\Mh,\zz)$ is generated by $\pi(A_i)$ and
$$H_1(\Mh,\zz)=\begin{cases}\zz_{4g+2}&\text{for $g$ even}\\
\zz_{8g+4}&\text{for $g$ odd}
\end{cases}$$
Now observe that since $B$ has order $2g+2$, $B^{g+1}$ is an involution. It is not
central in $\Mh$, so it is conjugate to $S$ or to $\ro S$. Therefore, using relations
\eqref{roro},\eqref{BB} and Lemma \ref{ConInv} we obtain
$$\pi(I_g)=\gen{\pi({\ro}),\pi({S}),\pi(\ro S)}=\gen{\pi({\ro}),\pi({B^{g+1}})}=\gen{2(2g+1),(g+1)(2g+1)}$$
Together with \eqref{eq:pi} this gives us desired result.
\end{dow}
\begin{tw}
The group $\Mhpm$ is generated by three symmetries.
\end{tw}
\begin{dow}
Let $\tau$ be the reflection across the $xz$-plane (Figure \ref{r1}) and $\eps_1=\tau B$,
$\eps_2=\tau A_{g+1}$. Since $\tau B\tau=B^{-1}$ and $\tau A_{g+1}\tau=A_{g+1}^{-1}$, both $\eps_1$
and $\eps_2$ are symmetries. Moreover $A_{g+1},B\in \gen{\tau,\eps_1,\eps_2}$, so
$\gen{\tau,\eps_1,\eps_2}=\Mhpm$.
\end{dow}



\end{document}